# A LAUNDRY SURFACE FOR BRAIDS

VICTOR A. NICHOLSON

*Department of Mathematical Sciences, Kent State University,*

*Kent, OH 44240, USA*

*victor@math.kent.edu*

ABSTRACT

Laundry surfaces for closed braid diagrams are presented. It is shown that braid diagrams are characterized by linking matrices obtained by lifting cycles from these surfaces. Oriented link types are then characterized by equivalence classes of linking matrices. Similar equivalence classes can be composed of Gordon and Litherland forms or Seifert matrices.

## 1. Introduction

A sequence can be formed by traveling around a knot diagram noting the crossing information and the order in which crossings are encountered. A laundry surface is about "order" in this sense. Combining a laundry surface with braid structure yields a family of surfaces such that both the surface and its embedding in 3-space are determined simply by the linking matrix. In Section 3 it is shown that oriented links are characterized by equivalence classes of linking matrices derived from such surfaces. This result is not surprising having been anticipated for many years. Most of the concepts are therefore familiar [4]. It is somewhat surprising that these matrices can represent braid diagrams. This is shown in Section 2. Laundry surfaces can be very twisted but are essentially "unknotted". The embeddings of laundry surfaces in 3-space are characterized by the laundry embedding theorem. This theorem is used in Section 4 to show that the laundry surface constructed from a braid diagram is ambient isotopic to its laundry form.

## 2. Matching Linking Matrices with Braid Diagrams

All spaces are assumed to be piecewise linear. A surface $S(L)$ constructed from a closed braid diagram $L$ is illustrated in Fig. 1. The diagram $L$ is a boundary component of $S(L)$ and in this example represents a Markov stabilzation move on a four-one knot. The surface construction begins with a plane annulus for each Seifert circle. The circle bands $C_1,...,C_4$ appear at the top of the annuli. Folds are made from left to right and twisted bands $X_1,...,X_5$ are attached at the braid crossings. Folds are also made on the right side and untwisted bands attached in a staircase pattern. The folds on the left and right are made in opposite directions. The surface is a disk-band surface in that it decomposes into a disk with bands attached. The bands added on the right side are considered to be part of the disk. The boundary of the surface has two components, the original diagram and an unknot that lies below. The removal of any circle band results in a surface whose boundary is the connected sum of these two components and is a link equivalent to the original diagram $L$.

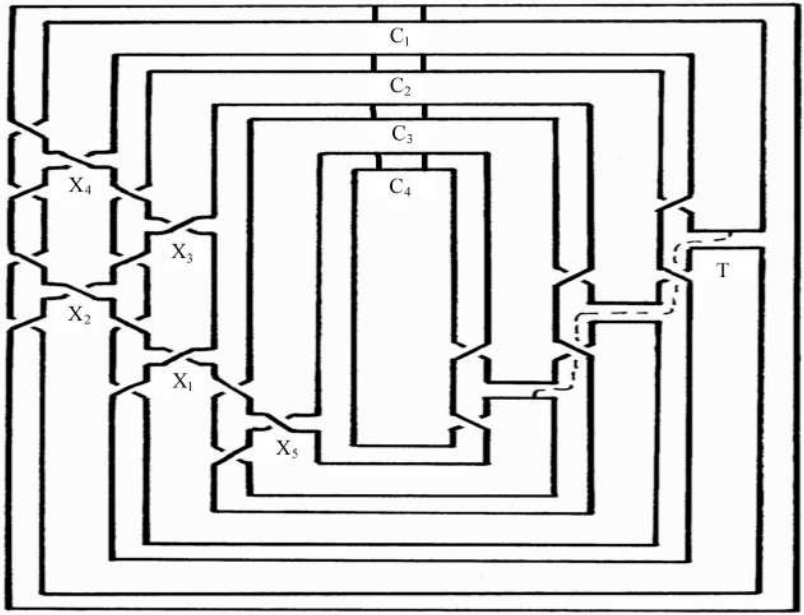

Fig. 1. Closed braid diagram with surface.

The surface $S'(L)$ illustrated in Fig. 2 (1) is said to be in laundry form (defined in Section 4). So called because the bands hang down like laundry when the surface is opened up. Theorem 3 in Section 4 shows that the surfaces $S(L)$ and $S'(L)$ are ambient isotopic. An elementary argument relying only on the braid structure would be to slide both the twisted bands and the circle bands down to the bottom of Fig. 1. Continue sliding and removing twists until only the twists in the twisted bands remain.

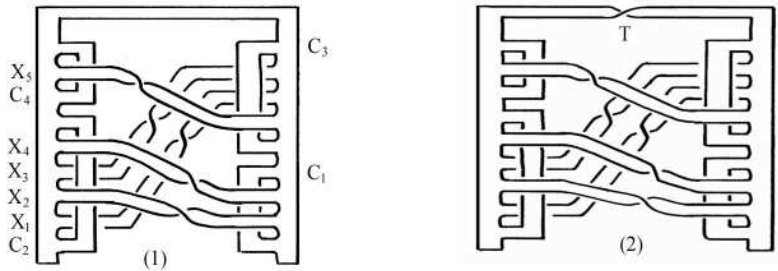

Fig. 2. Surfaces in laundry position.

The concentric Seifert circles are numbered from left to right. The columns between circles are also numbered from left to right. The surface $S'(L)$ in laundry form is similar to a three-page book with the circle bands on the middle page, the bands in odd columns above, and the bands in even columns below. The bands in a braid column are parallel to each other and this pattern is repeated in $S'(L)$. Each band is considered to have two feet (a first and a second foot) where it attaches to the disk. Laundry order is the order in which the bands appear around the outside of Fig. 2 (1). That is $\{C_2, X_1, ..., C_1\}$. In Fig. 1 a foot of each twisted band is encountered when traveling up the even circles. Let $f_i$ denote the list of the band feet in the order encountered traveling up the even circle $C_i$ and $u_i$ denote the foot of the circle band for the i-th circle. These labels define the first feet of the bands. Let $s_i$ denote the list of the band feet in the order encountered traveling down the odd circle $C_i$ and $v_i$ denote the foot of the circle band for the i-th circle. The laundry ordering of the band feet is given by an increasing list of the even circles interspersed with the lists of the first feet followed by a decreasing list of the odd circles interspersed with the lists of the second feet. That is, $(u_2 f_2 v_2 u_4 f_4 v_4 ...)(...u_3 s_3 v_3 u_1 s_1 v_1)$.

The first homology group of the surface is generated by a family of cycles, also denoted $\{C_2, X_1, ..., C_1\}$. Each cycle is oriented counter clockwise in Fig. 2. The linking matrix is $M(L)=(lk(c_i',c_j))$ where $c_i'$ is $c_i$ pushed off in "both directions", for $c_i, c_j \in \{C_2, X_1, ..., C_1\}$. It can be seen from Fig. 2 that the only nonzero entries in the linking matrix arise from one of the two properties. (L1) Suppose $X$ is the cycle for a twisted band in the i-th column. Then $lk(X', C_i) = -1$ and $lk(X', C_{i+1})=1$. (L2) Suppose $X$ and $Y$ are cycles for twisted bands which meet the same circle and the band for $X$ is above and left or below and right of the band for $Y$ then $lk(X',Y)=1$. The linking matrix $M(L)$ for the example representing the four-one knot $L$ appears in Fig. 3.

$$\begin{bmatrix} 0 & -1 & 1 & -1 & 1 & 0 & 0 & 0 & 0 \\ -1 & -1 & 1 & 0 & 1 & 0 & 1 & 1 & 0 \\ 1 & 1 & 1 & 0 & 0 & 0 & 0 & 0 & -1 \\ -1 & 0 & 0 & -1 & 1 & 0 & 1 & 1 & 0 \\ 1 & 1 & 0 & 1 & 1 & 0 & 0 & 0 & -1 \\ 0 & 0 & 0 & 0 & 0 & 0 & 1 & 0 & 0 \\ 0 & 1 & 0 & 1 & 0 & 1 & 1 & -1 & 0 \\ 0 & 1 & 0 & 1 & 0 & 0 & -1 & 0 & 0 \\ 0 & 0 & -1 & 0 & -1 & 0 & 0 & 0 & 0 \end{bmatrix}$$

Fig. 3. Linking matrix.

Two braids which are equivalent by the relation (B0) $\sigma_i \sigma_j = \sigma_j \sigma_i$ when $|i-j|>1$ are isotopic in the plane and have the same braid diagram. Let $D$ denote the set of braid diagrams.

**Theorem 1.** *The correspondence $L \to M(L)$ defines a bijection $M \to D$ of the set of linking matrices onto the set of braid diagrams.*

**Proof.** Begin with a linking matrix and recover the braid diagram as follows. The twist in a band is the diagonal entry. The zeros on the diagonal determine the number of circles and identify the circle bands in laundry order. Property L1 determines the column for each band. Consider the first column in the example matrix. The twisted bands are partitioned into sets all of which meet a common even circle by the lists of first feet. The bands within the list are encountered going up the even circle. Property L2 then determines the order of the bands on each odd circle. □

The removal of any row and column corresponding to a circle band from $M(L)$ results in a matrix that represents a Gordon and Litherland form [3] for the diagram. Let $F$ be the set of Gordon and Litherland forms obtained by deleting the last row and column for each matrix in $M$ together with the empty matrix. The last row and column correspond to the first circle. If $L$ consists of only a single circle then deleting the last row and column changes the one by one zero matrix to the empty matrix.

The laundry surface can be made orientable by adding a twisted band $T$ at the top as shown in Fig. 2 (2). In Fig. 1, cut the surface along the dotted line and insert a negative twisted band $T$. This results in an orientable surface as shown in Fig. 4.

Let $M'(L)=(lk(c_i^+,c_j))$ where $c_i^+$ is $c_i$ pushed off the positive side, for $c_i,c_j \in \{C_2, X_1, ..., C_1\}$. The positive side is the top on the right side of Fig. 2 (2). Let $M'$ denote the set of these matrices. The removal of any row and column corresponding to a circle band from $M'(L)$ results in a Seifert matrix for the diagram $L$. Let $S$ be the set of Seifert matrices obtained by deleting the last row and column for each matrix in $M'$ together with the empty matrix.

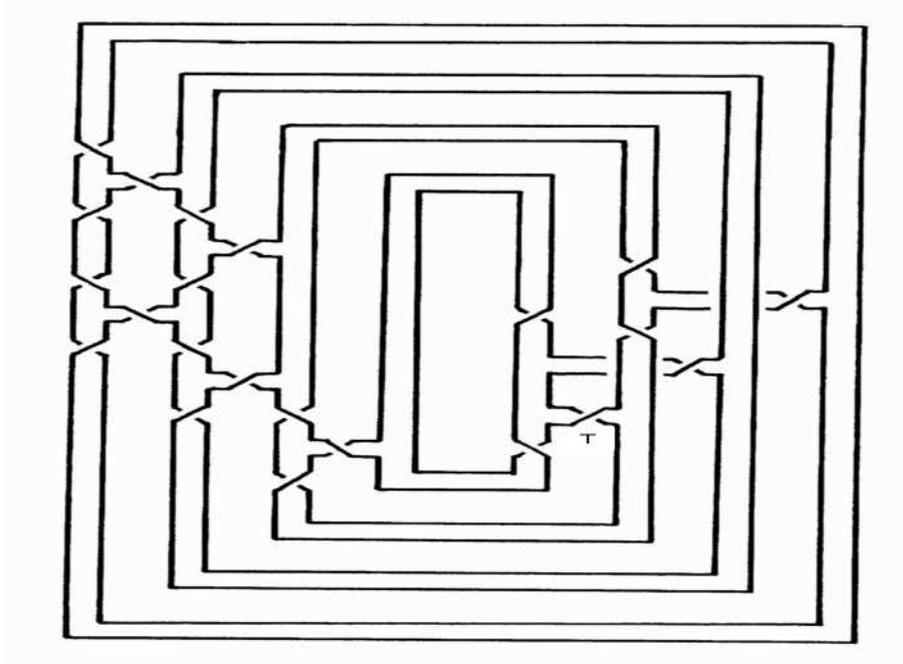

Fig. 4. Braid diagram with orientable surface.

**Corollary 1.** *There are bijections M′→D, F→D and S→D of the linking matrices for the orientable surfaces, the Gordon and Litherland forms, and the Seifert matrices, respectively, onto the set of braid diagrams.*

**Proof.** The cycles for the twisted bands pass through the band $T$ in the orientable surface. Define a matrix $N(L)=(n(c_i,c_j))$ where $n(c_i,c_j)=1$ if $c_i$ and $c_j$ are twisted bands and zero, otherwise, for $c_i,c_j \in \{C_2, X_1, ..., C_l\}$. The map $M′→M$ given by $M(L)=M′(L)+M′(L)^T+N(L)$ is a bijection because $M(L)$ determines $L$. It follows from condition (L1) that any row and column in $M(L)$ for a circle band is minus the sum of the rows and columns of the other circle bands. Thus $M$ can be recovered from $F$ so $M→F$ is a bijection. Let $S→F$ be given by $F(L)=S(L)+S(L)^T+N(L)$. The map $M′→S→F$ is the same as the bijection $M′→M→F$ so $M′→S$ is also a bijection. □

### 3. Matching Matrix Equivalence Classes with Oriented Links

Define two closed braid diagrams in $B$ to be equivalent if they are related by a finite sequence of braid diagrams such that each is obtained from its predecessor by one of the braid relations: (B1) the type II Reidemeister move; (B2) the Markov stabilization/destabilzation move; (B3) conjugacy; or (B4) the type III Reidemeister move. Let $B_E$ denote the set of equivalence classes of braid diagrams.

The corresponding matrix relations M1-M4 on the set $M$ of linking matrices are given below. These relations are a subset of the familiar relations known as *S*-equivalence. Two matices are *S*-equivalent if they are connected by a sequence of matrices such that each is obtained from its predecessor by either a single tube addition/deletion or a unimodular congruence. A unimodular congruence is a relation of the form $B=P^TAP$, where $detP=\pm 1$.

The moves are illustrated in laundry position in Figs. 5 and 6 where $C_i$ is an odd circle. There are four versions for each move. Versions for even circles can be obtained by rotating the figures 180 degrees in the plane of the paper and reversing all band crossings. Versions with different twists can be obtained by reversing all of the twists in a figure. The matrix relations are not affected by changing versions.

(M1) The type II Reidemeister move is illustrated in Fig. 5 (1). The move is addition/deletion of bands $X_1$ and $X_2$.

(M2) The Markov stabilization/destabilzation move is illustrated in Fig. 5 (2). The move is addition/deletion of bands $X_2$ and $C_{i+1}$. This move can also be obtained by tube addition/deletion (M1) accompanied by a slide of a foot of $X_1$ over $X_2$. The matrix $P$ then uses the cycle equation $C_{i+1}=X_1-X_2$ and a permutation interchanging $C_{i+1}$ and $X_2$.

(M3) The conjugation move is illustrated in Fig. 5 (3). The band $X_2$ can be obtained from the band $X_1$ by sliding the feet of $X_1$ over the bands $C_i$ and $C_{i+1}$. The matrix $P$ is the product of a unimodular matrix using the cycle equation $X_1-X_2=C_i+C_{i+1}$ and a permutation returning the matrix to laundry order.

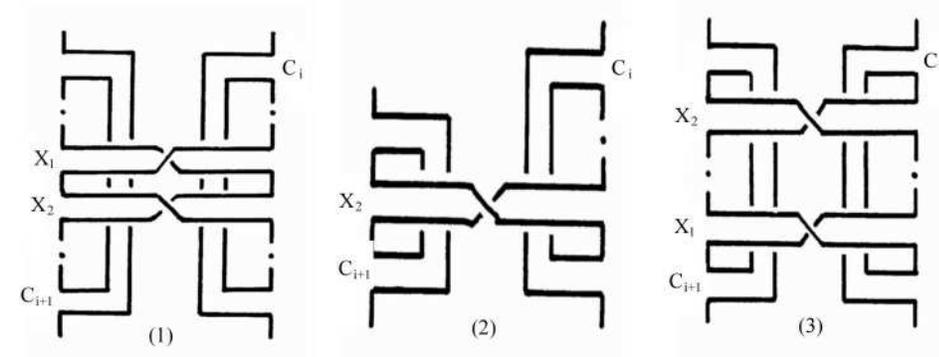

Fig. 5. Matrix moves M1, M2, and M3.

(M4) The type III Reidemeister move is illustrated in Fig. 6. The band $X_4$ can be obtained from the band $X_1$ by sliding the feet of band $X_1$ over the bands $X_2$ and $X_3$. The matrix $P$ uses the cycle equation $X_1+X_4=X_2+X_3$ and is the product of a unimodular matrix and a permutation returning the matrix to laundry order.

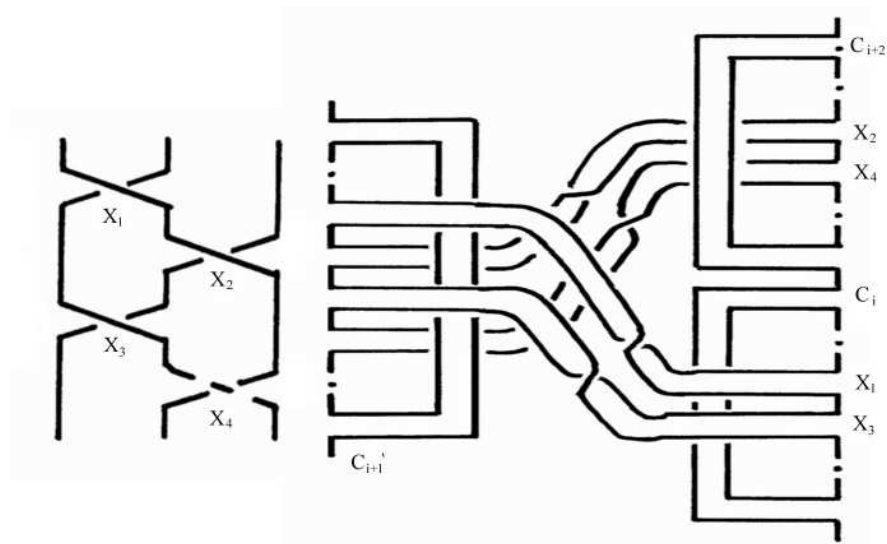

Fig. 6. Matrix move M4.

Define two matices in $M$ to be equivalent if they are related by a sequence of matrices in $M$ such that each is obtained from its predecessor by one of the moves M1-M4. Let $M_E$ denote the set of equivalence classes. Let $L_E$ denote the set of oriented link types in oriented 3-space.

**Theorem 2.** *There is a bijection $M_E \to L_E$ of the matrix equivalence classes onto the set of oriented link types.*

**Proof.** There is a surjection $B \to L_E$ of braid diagrams onto the set of oriented link types (Alexander's Theorem) [1,8]. The surjection $B \to L_E$ induces a bijection $B_E \to L_E$ (Markov's Theorem) [2,7]. If two matrices are related by one of the moves M1-M4 then the bijection $M \to B$ of Theorem 1 matches the matrices with equivalent braid diagrams. Conversely, if two braid diagrams are related by one of the braid moves then they are matched to equivalent matrices. Thus the map $M \to B$ induces a bijection $M_E \to B_E$. □

The set of equivalence classes $M_E$ is a complete invariant of oriented links.

**Corollary 2.** *There are bijections $M_E' \to L_E$, $F_E \to L_E$ and $S_E \to L_E$ of the matrix equivalence classes for the orientable surfaces, the Gordon and Litherland forms, and the Seifert matrices, respectively, onto the set of oriented link types.*

**Proof.** The matrix moves do not involve the twisted band $T$ used to make the laundry surface orientable. Thus $M_E$ can be replaced by $M_E'$. The set of Gordon and Litherland forms $F$ and the set of Seifert matrices $S$ match $M$ and $M'$, respectively. The equivalence relation can be defined on these sets by this matching. The cycle equation for a conjugation in column one changes to $X_1 - X_2 = -(C_2 + ...) + C_2$ because the row and column for the circle $C_1$ are deleted. The other matrix relations are unchanged. □

### 4. Laundry Surfaces

Fig. 7 shows the boundary the disk in Fig. 1 together with arcs through the bands. This is an example of a circle-with-chords. Consider a graph $G$ consisting of a circle subdivided by a set of points $\{a_1,...,a_n, b_1,...,b_n\}$ together with a set of $n$ oriented chords $E_i = a_i b_i$, $i=1,...,n$. The points are labeled so that $a_0$ and $b_0$ are adjacent on the circle, the oriented arc $J = a_0 b_0$ in the circle contains all of the vertices, $a_0 < a_1 < ... < a_n$, and $a_i < b_i$, for $i=0,...,n$. Let $E_0 = a_0 b_0$ (also called a chord) be the closure of the complement of $J$ in the circle. A circle-with-chords is a pair $(G,J)$. Two graphs $(G,J)$ and $(G',J')$ are said to be the same if there is a homeomorphism $h:(G,J) \to (G',J')$ that preserves the orientations of the arcs $J$ and $J'$. A circle-with-chords is essentially a based and oriented Gauss diagram. Start at the point $a_0$ in Fig. 7 and follow the arc $J$ around the circle to $b_0$. Note that the endpoints of the chords are encountered in the laundry order given in Section 2.

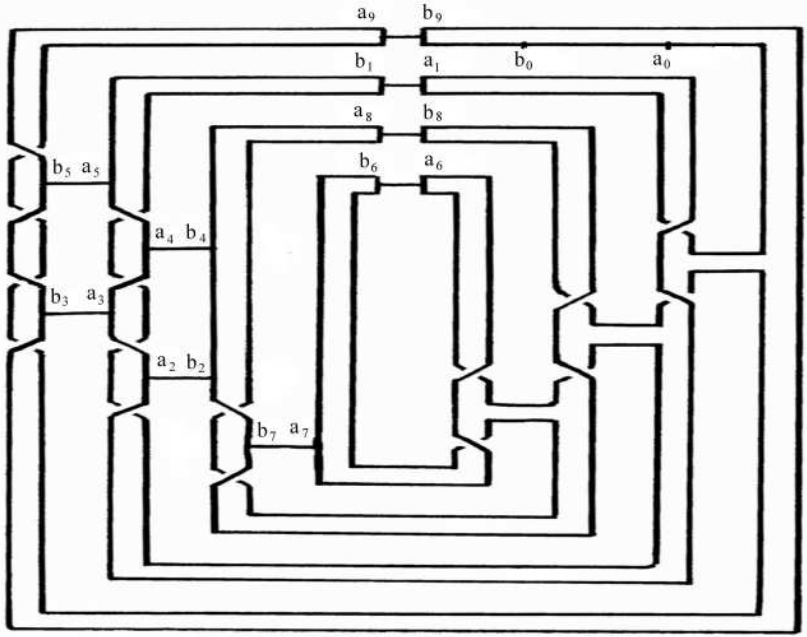

Fig. 7. Circle-with-chords.

Fig. 8 (1) shows a circle-with-chords that is the same as the one in Fig. 7. Suppose $G$ is in $R^3$, $p: R^3 \to R^2$ is the orthogonal projection, and $P_i=p(E_i)$, $i=0,...,n$, are the projected chords. The graph $(G,J)$ is said to be in laundry position if (1) every vertical line intersects $G$ in zero, one, or two points and, if two, neither of them is a vertex of $G$, (2) the arc $J$ of $G$ lies in the x-axis, (3) each $P_i$ is an arc in the bottom half of the xy-plane, (4) any pair of $P_i$ intersect in at most one point, (5) no $P_i$ meets the arc of $G$ at an interior point of $P_i$, and (6) if $P_i$ meets $P_j$ and $P_k$ with $i < j < k$, then $P_i$ meets $P_j$ first and then $P_k$. The central feature of laundry position is condition (6) which requires that the order in which a chord meets subsequent chords is the same as their order along the arc. Consider the face of a tennis racquet. The arc $J$ defines the "laundry line". The graph in Fig. 8 (1) is moved to laundry position by lifting the left and right sides to the x-axis.

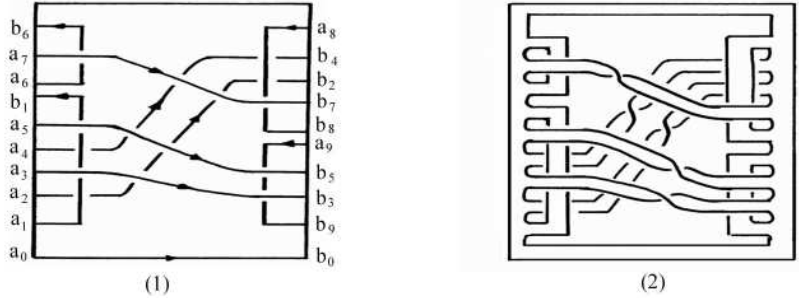

Fig. 8. Graph and laundry surface in laundry position.

A laundry surface is a triple $(S,G,J)$ consisting of a 2-manifold with boundary $S$ in $S^3$, a circle-with-chords $G$, and its arc $J$ satisfying the condition that $S$ is a regular neighborhood of $G$ in $S$ and is equivalent to a triple $(S',G',J')$ with $(G',J')$ in laundry position. Where triples $(S,G,J)$ and $(S',G',J')$ are considered to be equivalent if there is an orientation preserving homeomorphism $h$ of $S^3$ onto itself such that $h(S,G,J)=(S',G',J')$ and $h$ preserves the orientations of the arcs $J$ and $J'$. That is, equivalence is ambient isotopy.

A laundry surface can arise naturally from a disk-band surface. Cut a small hole in the disk, let the circle go around the hole and the chords pass through the bands. Recall that the constructed surface $S(L)$ in Fig. 1 is a disk-band surface. Let $S=S(L)-intB$, where $B$ is a small disk in the interior of its disk. Then $S$ is a regular neighborhood in $S$ of the circle-with-chords $(G,J)$ in Fig. 7 (adjusted slightly into the surface). It can be seen that $(G,J)$ is equivalent to the circle-with-chords in laundry position illustrated in Fig. 8 by sliding all of the chords down to the bottom of Fig. 7 and removing crossings. Thus $(S,G,J)$ is a laundry surface.

For $i=0,...,n$, let $I_i=a_ib_i$ a subarc of $J$. Let $X$ denote the family of cycles $X_i = E_i \cup I_i$, $i=0,...,n$, oriented in the direction of the chord $E_i$. The cycles in $X$ are a basis for the first homology group of $G$ and $S$. The linking matrix is $M=(lk(X_i',X_j))$, where $X_i'$ is $X_i$ pushed off in "both directions". Note that the oriented arc $J$ determines the cycle basis and the order of the rows and columns in the matrix. As noted in Theorem 1, the linking matrix was sufficient to determine the braid and the surface. However, in general, the linking matrix is not sufficient to determine equivalence of laundry surfaces. Consider the following example. Let $(G_1,J_1)$ be the graph in Fig. 8 (1) with the chords ($a_2b_2$, $a_3b_3$, $a_4b_4$, $a_5b_5$, and $a_7b_7$) for the twisted bands removed. Let $S_1$ be the surface in Fig. 8 (2) with the twisted bands removed. Then $(S_1,G_1,J_1)$ is a laundry surface. Let $S_2$ be surface obtained from $S_1$ by adding a negative twist to the top band (as in Fig. 2 (2)) and a positive twist to the bottom band. Then $(S_2,G_1,J_1)$ is also a laundry surface. The pairs $(S_1,G_1)$ and $(S_2,G_2)$ are ambient isotopic but $J_1$ cannot be carried to itself.

Suppose $e$ is an edge in $J$ both of whose vertices have degree three. Such as the edge $e=b_6a_8$ in the preceding examples. Suppose $D$ is a small regular neighborhood of $e$ in $S$. Place the label "a" on both points in $\partial D \cap J$ and the label "b" on both points where $\partial D$ meets the chords. The cyclic order of these points in $\partial D$ is either (aabb) or (abab). Define the turn at $e$ to be zero if the order is (aabb) and one, otherwise. The turn at an edge is an invariant of ambient isotopy for triples. It isn't necessary to know the turn at every edge. The arcs $I_i$ and $I_j$ are said to overlap if $I_i \cap I_j \neq \emptyset$,

$I_i \not\subset I_j$, and $I_j \not\subset I_i$, . Each pair of overlapping arcs with $1 \leq i < j \leq n$ defines an ordered pair $(i,j)$ called an overlapping pair. The overlap graph O for G has vertex set $\{1,...,n\}$ and edge set $\{ij \mid (i,j)$ is an overlapping pair $\}$. The set of interior first-edges of G is $\{uv \subset J \mid v=a_i$ for $i=min\{j \mid j \in C\}$, for each component C of O not containing one$\}$. The overlap graph for the above examples has four vertices and no edges making the edge $e$ an interior first-edge. Two surfaces $(S,G,J)$ and $(S',G',J')$ are said to have the same turns if they have the same turn at each interior first-edge.

A proof of the laundry embedding theorem and its corollary are given in [6].

**Theorem** (laundry embedding) *Two laundry surfaces are equivalent iff they have the same graph, linking matrix, and turns.*

**Corollary** *Suppose N and N' are 2-manifolds with boundary (possibly empty). Suppose $\{B_1,...,B_r\}$ and $\{B_1',...,B_r'\}$ are families of pairwise disjoint disks in N and N', respectively. Suppose $N-\cup intB_i=(S,G,J)$ and $N'-\cup intB_i'=(S',G',J')$ are laundry surfaces. If $h:(S,G,J) \simeq (S',G',J')$ is an equivalence such that $h_1(\partial B_i)=\partial B_i'$, for $i=1,...,n$ then $h_1$ can be assumed to carry N to N'.*

Remark. The doubled-delta move [5] is illustrated in Fig. 9. The move preserves the graph, linking matrix, and turns of a surface. The laundry embedding theorem implies that a doubled-delta move must result in either an equivalent surface or one that is not a laundry surface. The bands can be arbitrarily numbered and a direction arbitrarily chosen for one of the bands. Note that the order in which this band meets the other two changes. If one figure satisfies condition (6) for laundry position then the other may not.

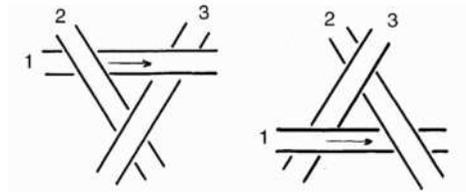

Fig. 9. Doubled-delta move.

**Theorem 3.** *Suppose L is a closed braid diagram. The surface $S(L)$ is ambient isotopic to the surface $S'(L)$ in laundry form.*

**Proof.** We use the corollary. As noted above, the surface $(S,G,J)$ in Fig. 1 is a laundry surface where $S=S(L)-intB$ and G is the graph in Fig. 7. The surface $(S'(L),G',J')$ in Fig. 8 (2) is also a laundry surface and $G'$ is the graph in Fig. 8 (1). The cycle bases $\{C_2, X_1, ...,C_1\}$ and the linking matrices are defined by the graphs. In Fig. 8 (1), the cycle $X_0$ goes around the outside rectangle. In Fig. 7, the cycle $X_0$ goes around the boundary of $B$. Both linking matrices have a row and column consisting entirely of zeros for the cycle $X_0$. Showing that the cycles $\{C_2, X_1, ...,C_1\}$ have the properties L1 and L2 on both surfaces will show that the two surfaces have the same linking matrix. The properties for $S'(L)$ were noted in the proof of Theorem 1. The cycles in $G'$ can be simplified by considering $B$ to be a small disk near the chord $a_0b_0$ in Fig. 7. The cycles in G are equivalent in S to cycles with the following properties: each cycle in $\{C_2, X_1, ...,C_1\}$ passes up the even circle, through the band and then back down the odd circle; the cycles for the circle bands are oriented clockwise for odd circles and counterclockwise for even circles. Contributions to linking arise only when both cycles pass through a twist in the surface. If both curves travel through the twist in the same direction the contribution is the sign of the twist. If they pass in opposite directions the sign changes. Properties L1 and L2 are established by inspection of the twists produced by the folding construction of the surface. Most of the twists cancel in pairs. This establishes that the linking matrices are the same. The turn is always zero for the surface $S=S(L)-intB$ because the cycle $X_0$ simply goes around the hole and the chords all point away from this hole. The turn is zero in Fig. 8 (2) because there are no twists on the laundry line. The laundry embedding theorem implies that these laundry surfaces are equivalent. The outside rectangle in Fig. 8 (2) bounds a disk behind the surface. The corollary then implies that the isotopy extends carrying the disk $B$ to this disk. Attaching a disk to the outside rectangle is equivalent to removing the band for $X_0$ and produces Fig. 2 (1). Thus $S(L)$ is ambient isotopic to the surface $S'(L)$. □

Remark. The orientable surfaces in Fig. 2 (2) and Fig. 4 can be shown to be equivalent using the proof of Theorem 3 and the replacement lemma for laundry surfaces in [6].